\documentclass[letterpaper, 10 pt, conference]{ieeeconf}  

\IEEEoverridecommandlockouts                

\usepackage{xcolor}

\usepackage{hyperref}
\usepackage{mathtools}          
\usepackage{amssymb}
\usepackage{graphicx}
\usepackage[ruled,vlined]{algorithm2e}
\usepackage{cuted}
\usepackage{pifont}

\newtheorem{lemma}{Lemma}[section]

\newtheorem{remark}{Remark}

\newcommand{\diff}{\,{\rm d}}

\newcommand{\R}{\,{\mathcal{R}}}

\newcommand{\Cay}{\,{\rm Cay}}
\newcommand{\Om}{\Omega}

\newcommand{\cmark}{\ding{51}} 
\newcommand{\xmark}{\ding{55}} 

\graphicspath{{figures/}}

\title{\LARGE \bf
	Cayley Commutator-free Methods for Krotov-Type Algorithms in Quantum Optimal Control
}

\author{Boris Wembe$^{1}$, Usman Ali$^2$, Torsten Meier$^3$ and Sina Ober-Blöbaum$^{4}$ 
	\thanks{*This work was not supported by any organization}
	\thanks{$^{1}$ Faculty of Electrical Engineering, Mathematics and Computer Science, Paderborn University, Paderborn, Germany
		{\tt\small wboris@math.upb.de}}%
	\thanks{$^{2}$ Institut für Physik, Carl von Ossietzky Universität, D-26111 Oldenburg, Germany
		{\tt\small usman.ali@uol.de}}%
	\thanks{$^{3}$ Institute for Photonic Quantum Systems (PhoQS) and Department of Physics, Paderborn University, Paderborn, Germany
		{\tt\small tmeier@uni-paderborn.de}}%
	\thanks{$^{4}$ Faculty of Electrical Engineering, Mathematics and Computer Science, Paderborn University, Paderborn, Germany
		{\tt\small sinaob@math.upb.de}}%
}

\begin{document}

	\maketitle
	\thispagestyle{empty}
	\pagestyle{empty}
	
	\begin{abstract}
		This paper presents a class of structure-preserving numerical methods for quantum optimal control problems, based on commutator-free Cayley integrators. Starting from the Krotov framework, we reformulate the forward and backward propagation steps using Cayley-type schemes that preserve unitarity and symmetry at the discrete level. This approach eliminates the need for matrix exponentials and commutators, 
		leading to significant computational savings while maintaining higher-order accuracy.  
		We first recall the standard linear setting and then extend the formulation to nonlinear Schrödinger and Gross--Pitaevskii equations using a Cayley--polynomial interpolation strategy. Numerical experiments on state-transfer problems illustrate that the CF--Cayley method achieves the same accuracy as high-order exponential or Cayley--Magnus schemes at substantially lower cost, especially for longtime or highly oscillatory dynamics.  
		In the nonlinear regime, the structure-preserving properties of the method ensure stability and norm conservation, making it a robust tool for large-scale quantum control simulations. The proposed framework thus bridges geometric integration and optimal control, 
		offering an efficient and reliable alternative to existing exponential-based propagators.
	\end{abstract}

	\section{Introduction}
	
	Quantum technologies are rapidly advancing in areas such as quantum computation, quantum communication, and quantum simulation. A central requirement for these applications is the \emph{ability to prepare and manipulate quantum states with high fidelity} under realistic physical constraints \cite{Kobylarov25}. This task is naturally formulated as a \emph{quantum optimal control problem} (QOCP), where one seeks external control fields that steer a quantum system, governed by the (linear or nonlinear) Schrödinger equation, towards a desired target state or unitary operation \cite{Glaser2015, Koch2022}. QOCPs were first studied to steer chemical reactions with tailored laser pulses, enabling bond-selective dissociation \cite{Brif2010}. These optimal-control techniques have since matured into general tools for quantum engineering. More recently, they have been applied to sculpt the phase-space distributions of Bose–Einstein condensates and to coherently transfer atoms into selective-Bloch bands of periodic-optical lattices generated by lasers fields \cite{Dupont2021, Dupont2023}.
	
	Among the various approaches to QOCP, \emph{gradient-based iterative methods}, in particular Krotov’s method have proven to be powerful and versatile. Originating from the work of Krotov in the mid-1990s \cite{Krotov1996}, the method has been systematically developed for applications in quantum control \cite{Grond2009, Hohenester2014, Koch2022}. It offers guaranteed monotonic convergence under suitable conditions and has been successfully applied in diverse contexts such as quantum state preparation, gate implementation, and nonlinear many-body dynamics \cite{Goerz2019, Jager2014, Reich2012}. However, its \emph{computational bottleneck lies in the repeated forward and backward propagation} of the Schrödinger equation across fine temporal grids. Since each iteration requires a full solution of the time-dependent Schrödinger dynamics, the efficiency of the underlying time integrator directly determines the scalability of the overall algorithm.  
	
	Classical integrators, such as high-order Runge–Kutta schemes or Magnus-type exponential integrators \cite{Blanes2009}, are widely used for this task. While accurate, they often incur substantial computational costs due to either the evaluation of matrix exponentials or the computation of nested commutators. For large-scale systems or highly oscillatory dynamics, these costs can become prohibitive, thereby limiting the applicability of optimal control to more complex quantum systems.  
	
	Recently, \emph{commutator-free Cayley methods} have been proposed as a new class of structure-preserving integrators for differential equations evolving on Lie groups \cite{Wembe2025}. Instead of relying on expensive exponential maps or commutator expansions, these schemes are built from compositions of \emph{Cayley transforms}, which provide a rational, unitarity-preserving approximation of the matrix exponential. As a result, they combine three essential features: \emph{exact unitarity of the propagator}, which is crucial for quantum dynamics; \emph{computational efficiency}, avoiding matrix exponentials and commutator evaluations; \emph{High accuracy}, 4th-order convergence achieved with modest effort.  
	
	These properties make commutator-free Cayley integrators \emph{natural candidates for accelerating quantum optimal control algorithms}. In particular, their efficiency and stability suggest that they can significantly reduce the cost of Krotov-type algorithms, especially in scenarios involving many iterations, fine temporal resolution, or nonlinear extensions of the Schrödinger equation.  
	
	In this work, we explore the integration of commutator-free Cayley methods into Krotov’s framework for quantum optimal control. Our focus is on \emph{quantum state preparation problems}, both in linear (few-level systems) and nonlinear (mean-field Gross–Pitaevskii-type) settings. We demonstrate that the Cayley-based approach achieves a favorable balance between accuracy, computational cost, and structural fidelity, thereby providing a promising pathway towards more scalable quantum control strategies.

	\section{Problem Formulation and Necessary Conditions}

	\subsection{Quantum Optimal Control Problem}
	\label{sec:QOCP}
	
	Quantum state preparation lies at the heart of quantum theory, since the outcome of any experiment or quantum information protocol depends critically on the ability to initialize well-defined states. In practice, however, preparing such states is experimentally challenging due to decoherence, unwanted couplings, and stringent constraints on the available control fields. Optimal control theory provides a systematic framework to overcome these obstacles by designing external fields that reliably steer the system toward the desired target state with high fidelity while respecting physical limitations.  
	For the purpose of this article, we consider a controlled quantum system evolving in a Hilbert space $\mathcal{H}\cong \mathbb{C}^d$, governed by the nonlinear Schrödinger equation
	\begin{align}
		& i \frac{\mathrm{d}}{\mathrm{d}t}\psi(t) 
		= \Big(H_0 + \sum_{j=1}^m u_j(t) H_j\Big)\psi(t) 
		+ F(\psi(t)), \nonumber \\
		&\hspace{6cm} \psi(0) = \psi_0.
		\label{eq:system_dynamics}
	\end{align}
	Here $H_0$ is the drift Hamiltonian, $H_j$ are control Hamiltonians driven by real control fields $u_j(t)\in L^2([0,T])$, and $F$ denotes a possible nonlinear term (e.g.\ the cubic interaction in the Gross–Pitaevskii equation). The state $\psi(t)$ is normalized, $\|\psi(t)\|=1$, reflecting probability conservation. The optimal control problem consists in minimizing the functional\footnote{Alternative formulations of the cost functional include a regularization term penalizing the time derivative of the control, 
		$J[u] = J_T(\psi(T)) + \tfrac{\gamma}{2}\int_0^T |\dot{u}(t)|^2\,\mathrm{d}t$. 
		Such terms promote smooth and physically feasible controls and can naturally enforce positivity when the control represents a nonnegative quantity, such as a laser intensity. A brief discussion of this approach is provided in Appendix~\ref{app:alt-cost}.
	}

	{\small
		\begin{equation}
			J[u] = \frac{1}{2}\left(1 - |\langle \psi_T, \psi(T)\rangle|^2 \right) + \sum_{j=1}^m \frac{\alpha_j}{2}
			\int_0^T |u_j(t)|^2 \,\mathrm{d}t,
		\end{equation}
	}
	where $\psi_T$ is the desired target state. The first term measures the fidelity between $\psi(T)$ and $\psi_T$, and the second term penalizes large control amplitudes with weights $\alpha_j>0$.  
	This setting is representative of state-to-state transfer tasks in quantum control, and can be generalized to gate optimization, ensemble control, or dissipative dynamics by adapting $\psi(t)$, $F$, and the cost functional.

	\subsection{Pontryagin Framework and Necessary Conditions}
	\label{sec:PMP}
	
	To characterize optimal solutions, according to the Pontryagin’s Maximum Principle \cite{Pontryagin1962, Krotov1996}, we introduce the adjoint state $\lambda(t)\in \mathbb{C}^d$ and the Hamiltonian functional
	\begin{equation}
		\mathcal{H}(\psi,\lambda,u) = 
		\text{Im} \,\langle \lambda, (H_0 + \sum_{j=1}^m u_j H_j)\psi + F(\psi) \rangle - g(u),
	\end{equation}
	where $g(u)=\sum_{j=1}^m \tfrac{\alpha_j}{2}|u_j|^2$. The necessary conditions for optimality are given by the two-point boundary value problem
	{\small
		\begin{subequations}
			\begin{align}
				i \dot{\psi}(t) &= \Big(H_0 + \sum_{j=1}^m u_j(t) H_j\Big)\psi(t) + F(\psi(t)), \nonumber \\
				&\hspace{5cm} \psi(0) = \psi_0,   \\
				i \dot{\lambda}(t) &= \Big(H_0 + \sum_{j=1}^m u_j(t) H_j\Big)\lambda(t) + \nabla_\psi F(\psi(t))^\ast \lambda(t),  \nonumber \\
				&\hspace{3.5cm} \lambda(T) = -\nabla_\psi \Phi(\psi(T)), \\
				\alpha_j u_j &= \text{Im} \langle \lambda(t), H_j \psi(t)\rangle. 
			\end{align}
		\end{subequations}
	}
	Thus, the optimality system couples forward propagation of the state $\psi$ with backward propagation of the adjoint $\lambda$, linked through pointwise stationarity conditions for the controls. Solving this coupled system is computationally demanding, motivating iterative schemes such as Krotov’s method.

	\subsection{Krotov-Type Iterative Method}
	\label{sec:Krotov}
	
	Krotov’s method provides an efficient iterative scheme 
	to solve the optimality system while guaranteeing monotonic decrease of the cost functional. Starting from an initial guess $u^{(0)}$, each iteration $k \mapsto k+1$ consists of three main steps, outlined in Algorithm~\ref{alg:Krotov}.  
	\begin{algorithm}[h!]
		\caption{Krotov Iteration Algorithm}
		\label{alg:Krotov}
		\KwIn{Initial guess $u^{(0)}(t)$, penalty parameters $\alpha_j>0$, target state $\psi_T$}
		\KwOut{Optimized control $u^{(k)}(t)$}
		
		\For{$k = 0,1,2,\dots$ until convergence}{
			\textbf{Forward propagation:} Solve the state equation with $u^{(k)}$ to obtain $\psi^{(k)}(t)$. \\
			\textbf{Backward propagation:} Solve the adjoint equation backward in time with $u^{(k)}$ to obtain $\lambda^{(k)}(t)$. \\
			\textbf{Control update:} Update each control component $u_j^{(k)}$ pointwise in time via
			\[
			u_j^{(k+1)}(t) 
			= u_j^{(k)}(t) 
			+ \frac{1}{\alpha_j} 
			\text{Im} \langle \lambda^{(k)}(t), H_j \psi^{(k+1)}(t)\rangle .
			\]
		}
	\end{algorithm}
	
	A detailed presentation of Krotov-type optimal control algorithms, including their derivation and monotonic convergence properties, can be found in \cite{Krotov1996, Reich2012, Goerz2019}. The algorithm guarantees $J[u^{(k+1)}] \le J[u^{(k)}]$ under mild assumptions, which makes it particularly robust in quantum applications. However, its numerical performance is strongly influenced by the accuracy of the forward and backward propagations. Standard explicit schemes may accumulate phase errors or fail to preserve unitarity, especially for long-time or highly oscillatory dynamics. This motivates the use of structure-preserving integrators, such as commutator-free Cayley methods, which retain the geometric properties of quantum evolution while reducing computational cost. In the next section, we introduce these integrators and analyze their role in enhancing 
	Krotov-type algorithms for quantum optimal control.

	\section{Cayley commutator-free method (CFCT)}
	
	\subsection{Commutator-Free Cayley Schemes in the Linear Case}
	\label{sec:CFC-linear}
	
	More details about the results presented in this section can be found in \cite{Wembe2025}, 
	while general background on Cayley methods is given in \cite{DLP:1998,IMNZ:2000,Iserles:2001}. See also the appendix.  
	Consider the time-dependent Schrödinger equation
	\begin{equation}
		\dot{\psi}(t) = A(t)\psi(t), 
		\quad A(t)^\dagger = -A(t) = -i H(t),
	\end{equation}
	whose exact solution is given by the unitary propagator
	\[
	\psi(t+\delta t) 
	= \mathcal{T}\exp\!\left(\int_t^{t+\delta t} A(s)\,\mathrm{d}s\right)\psi(t),
	\]
	with $\mathcal{T}$ the time-ordering operator.  
	Direct computation of the exponential is not always possible and expensive, while standard polynomial integrators may fail to preserve unitarity.  
	To address this, Magnus-type expansions approximate the solution as the exponential or Cayley transform of a series $\Omega(t) = \sum_{m=0}^\infty \Omega_m$; see, for instance, \cite{Iserles:2001}, where fourth- and sixth-order Cayley methods were developed. However, these approaches generally require nested commutators, which become costly in high dimensions. To avoid this, commutator-free methods have been proposed \cite{AF:2011,BM2006}, achieving high-order accuracy without commutators, though originally in the exponential setting.  
	Recently, we combined Cayley–Magnus expansions with the commutator-free idea to construct a fourth-order \emph{commutator-free Cayley (CFC)} method \cite{Wembe2025}. This scheme retains the structure-preserving advantages of Cayley integrators, but at significantly reduced computational cost. In Section~\ref{sec:CFC-Krotov} we illustrate the performance gain obtained by incorporating this method into Krotov’s algorithm for quantum optimal control.  
	
	For a given operator $A$, the Cayley transform is defined as
	\begin{equation}
		\Cay(A(\delta t)) 
		= \left(I + \tfrac{1}{2} A(\delta t)\right)^{-1}
		\left(I - \tfrac{1}{2} A(\delta t)\right),
	\end{equation}
	which is unitary, symmetric, and coincides with the Crank–Nicolson scheme.  
	The CFC idea is to approximate the propagator by a product of Cayley transforms,
	\begin{equation}
		U_{\mathrm{CF\text{-}Cay}}(\delta t, A) 
		= \prod_{\ell=1}^s \Cay(\delta t, \widetilde{A}_\ell),
	\end{equation}
	where the effective operators $\widetilde{A}_\ell$ are linear combinations of $A$ evaluated at quadrature nodes in $[t,t+\delta t]$.  
	As an example, a symmetric three-stage ($s=3$) fourth-order CFC scheme reads
	\begin{equation}
		U^{[4]}_{\mathrm{CFC}}(\delta t) 
		= \Cay(\delta t, \widetilde{A}_1)\,
		\Cay(\delta t, \widetilde{A}_2)\,
		\Cay(\delta t, \widetilde{A}_3),
	\end{equation}
	with
	\begin{align}
		\widetilde{A}_1 &= \alpha_{11} A(t+c_1 \delta t) 
		+ \alpha_{12} A(t+c_2 \delta t), \nonumber \\
		\widetilde{A}_2 &= \alpha_{21} A(t+c_1 \delta t) 
		+ \alpha_{22} A(t+c_2 \delta t), \\
		\widetilde{A}_3 &= \alpha_{31} A(t+c_1 \delta t) 
		+ \alpha_{32} A(t+c_2 \delta t), \nonumber
	\end{align}
	where the nodes and coefficients are
	{\small
		\[
		\begin{cases}
			c_{1,2} = \tfrac{1}{2} \mp \tfrac{\sqrt{3}}{6}, 
			~~ \alpha_{22} = 0, \\[0.3em]
			\alpha_{11} = \dfrac{2^{1/3}}{3} + \dfrac{2^{2/3}}{6} + \dfrac{2}{3},
		\end{cases} ~
		\begin{cases}
			\alpha_{21} = 1 - 2\alpha_{11}, \\
			\alpha_{12} = \alpha_{11} - \alpha_{11}^2, \\
			\alpha_{31} = \alpha_{11}, ~ \alpha_{32} = - \alpha_{12}. 
		\end{cases}
		\]
	}
	This scheme achieves fourth-order accuracy while each factor remains a Cayley transform, and is therefore unitary by construction. Compared to Magnus- or exponential-based commutator-free methods, the CFC scheme avoids commutator growth and expensive exponentials, making it particularly efficient for repeated forward/backward propagations in optimal control algorithms.

	\subsection{Commutator-Free Cayley Schemes in the Nonlinear Case}
	\label{sec:CFC-nonlinear}
	
	We now extend the commutator-free Cayley approach to nonlinear problems of the form
	\begin{equation}
		Y'(t) = A(t, Y(t))\,Y(t), 
		\qquad Y(t_0) = Y_0 \in G,
	\end{equation}
	where $G$ is a quadratic Lie group with associated Lie algebra~$\mathfrak{g}$.  
	Here the operator $A(t,Y)$ depends not only on time but also on the current solution, 
	as in nonlinear Schrödinger or Gross--Pitaevski-type equations. Our goal is to extend the CF-Cayley method in this context, combining it with polynomial interpolation, while preserving the order of approximation.  
	
	\paragraph{Algorithmic idea}
	Given a collection of previously computed approximations 
	$(t_0,Y_0),\dots,(t_{k-1},Y_{k-1})$, 
	we use Lagrange interpolation to approximate the trajectory $Y(t)$ 
	within the most recent $k$ steps.  
	Let $p_{k-1}$ denote the interpolating polynomial such that
	\[
	p_{k-1}(t_i) = Y_i, \qquad i=0,\dots,k-1.
	\]
	Evaluating this polynomial at two Gauss--Legendre nodes 
	$x_1 = \tfrac{1}{2} - \tfrac{\sqrt{3}}{6}$ and 
	$x_2 = \tfrac{1}{2} + \tfrac{\sqrt{3}}{6}$ 
	yields the stage values
	\[
	\begin{aligned}
		p_{k-1}(t_{k-1}+hx_1) \approx Y(t_{k-1}+hx_1), \\
		p_{k-1}(t_{k-1}+hx_2) \approx Y(t_{k-1}+hx_2),
	\end{aligned}
	\]
	which serve as inputs for the nonlinear evaluation of $A(t,Y)$ within each Cayley step.
	
	\paragraph{One-step update}
	The next approximation $Y_k$ is computed using a composition 
	of Cayley transforms that mirrors the linear commutator-free scheme:
	\begin{equation}
		\begin{aligned}
			Y_k = \Cay\!\left(\alpha_{11}\widetilde A_1(h) 
			+ \alpha_{12}\widetilde A_2(h)\right)
			\Cay\!\left(\alpha_{21}\widetilde A_1(h)\right) \\
			\Cay\!\left(\alpha_{11}\widetilde A_1(h)
			- \alpha_{12}\widetilde A_2(h)\right) Y_{k-1},
		\end{aligned}
	\end{equation}
	where
	\[
	\begin{aligned}
		\widetilde A_1(h) = A\!\left(t_{k-1}+hx_1, p_{k-1}(t_{k-1}+hx_1)\right), \\
		\widetilde A_2(h) = A\!\left(t_{k-1}+hx_2, p_{k-1}(t_{k-1}+hx_2)\right),
	\end{aligned}
	\]
	and the coefficients $\alpha_{ij}$ are identical to those used in the linear $4th$-order CF--Cayley scheme.
	
	After each step, the interpolation polynomial is updated to include the newly computed value $(t_k,Y_k)$, thus enabling continuous reuse of the last $k$ states for subsequent evaluations.  
	This procedure yields a symmetric, unitary integrator that preserves the Lie group structure even in the presence of nonlinear dependencies. 
	
	\paragraph{Cayley--Polynomial Algorithm.}
	The complete method, referred to as \emph{CaylPol}, 
	is summarized in Algorithm~\ref{alg:CaylPol} below.
	
	\begin{algorithm}[h!]
		\caption{CaylPol Integrator with $k$-step startup}
		\label{alg:CaylPol}
		\KwIn{Initial time $t_0$, step size $h$, number of steps $M$, startup size $k\!\ge\!2$, initial state $Y_0$} \KwOut{Approximations $\{Y_i\}_{i=0}^M$}
		
		\tcp{Startup phase: compute first $k$ states with a standard scheme}
		$\text{approx} \leftarrow \text{Startup}(k,t_0,h,Y_0)$\;
		$Y \leftarrow \text{approx}[k-1]$\;
		\tcp{Initialize interpolation data}
		$\text{times} \leftarrow [t_0,\dots,t_{k-1}]$\;
		$\text{states} \leftarrow [\text{approx}[0],\dots,\text{approx}[k-1]]$\;
		
		\For{$n=0$ \KwTo $M-k$}{
			$t^\star \leftarrow t_{k-1+n}$\;
			$p_1 \leftarrow \text{Interpolation}(t^\star + hx_1, \text{times}, \text{states})$\;
			$p_2 \leftarrow \text{Interpolation}(t^\star + hx_2, \text{times}, \text{states})$\;
			$Y \leftarrow \text{CayleyFreeStep}(t^\star,Y,p_1,p_2,h)$\;
			append $Y$ to $\text{approx}$\;
			update $(\text{times},\text{states})$ with the latest $k$ values\;
		}
		\Return $\text{approx}$\;
	\end{algorithm}
	
	\begin{remark}
		The CaylPol algorithm generalizes commutator-free Cayley schemes to nonlinear systems on Lie groups. By combining polynomial interpolation and Cayley composition, it maintains both unitarity and higher accuracy without requiring nested commutators. This makes it a suitable candidate for nonlinear quantum control problems, where the system dynamics often depend on the state amplitude itself. Notice that the startup of this algorithm requires the computation of some few first term with a standard scheme such as Runge-Kutta-Munthe-Kaas.
	\end{remark}

	\section{CF-Cayley Methods in Quantum Optimal Control}
	\label{sec:CFC-Krotov}
	
	\subsection{Incorporation into Krotov’s Algorithm}
	\label{sec:Krotov-CFC}
	
	The Krotov method requires repeated forward propagation of the state $\psi^{(k)}(t)$ 
	and backward propagation of the adjoint $\lambda^{(k)}(t)$ at each iteration.  
	Replacing the standard integrator with a commutator-free Cayley (CFC) scheme ensures 
	that each propagation step is unitary (or norm-preserving in the nonlinear case) and symmetric, 
	thereby reducing the accumulation of numerical errors across iterations.

	\begin{algorithm}[h!]
		\caption{Krotov Iteration with CF-Cayley Propagation in linear case}
		\label{alg:Krotov-CFC}
		\KwIn{Current control $u^{(k)}(t)$, state $\psi_0$, target $\psi_T$, penalty parameters $\alpha_j>0$}
		\KwOut{Updated control $u^{(k+1)}(t)$}
		
		\For{$k = 0,1,2,\dots$ until convergence}{
			\textbf{Forward propagation:} Propagate $\psi^{(k)}$ using a CF-Cayley method. \\
			\textbf{Backward propagation:} Propagate $\lambda^{(k)}$ backward in time using the same CF-Cayley scheme. \\
			\textbf{Control update:} Update each $u_j^{(k)}$ pointwise via
			\[
			u_j^{(k+1)}(t) = u_j^{(k)}(t) + \frac{1}{\alpha_j} \text{Im} \langle \lambda^{(k)}(t), H_j \psi^{(k+1)}(t)\rangle .
			\]
		}
	\end{algorithm}
	
	This modification leaves the monotonic convergence property of Krotov’s method intact, 
	but improves numerical stability and accuracy, especially for long-time or highly oscillatory dynamics.

	\begin{remark}
		In the nonlinear case, the adjoint equation depends explicitly on the forward state trajectory. Employing a higher-order integrator for the backward propagation would therefore require either storing additional intermediate states during the forward propagation or recomputing them on demand, both of which increase computational cost. An alternative is to use a lower-order but structure-preserving scheme for the adjoint propagation, which often provides a good balance between accuracy and efficiency. A detailed comparison of these approaches lies beyond the scope of the present work.
	\end{remark}

	\subsection{Example 1: Linear Schrödinger Equation for Non-interacting Trapped Cold Atoms in a Driven Optical Lattice}
	\label{sec:linear_example}
	
	As a first benchmark, we consider the case of a linear Schrödinger equation of the form
	\begin{equation}
		\begin{aligned}
			i\partial_t \psi(x,t) &= \left(-\Delta + V(x) + u(t) W(x)\right)\psi(x,t), \\
			&\hspace{3.4cm} \psi(x,0)=\psi_0(x),
		\end{aligned}
	\end{equation}
	%
	
	where the combined potential $V(x) =  V_0\sin^2\left(\frac{\pi x}{d}\right) + \frac{m\omega_T^2 x^2}{2}$ represents the superposition of a one-dimensional optical lattice and a quadratic trapping potential. The periodic optical lattice potential induces tunneling between modes, while the harmonic potential provides an overall confinement. The driving term $u(t)W(x)$ accounts for the time-dependent modulation of the control operator $W(x) = m \omega_T^2 x^2/2$ by the control field $u(t)$. This setup naturally incorporates key elements commonly realized in cold-atom experiments, where optical lattices and harmonic traps are often combined to manipulate atomic motion with high precision. The interplay between the optical lattice, the trapping potential, and the time-periodic driving field enables studies of coherent quantum transport. In particular, such a configuration has been explored and proposed as a platform to realize various transport dynamics associated with localized modes of the system \cite{Cao2020,Ali2023,Ali2024}.  
	We demonstrate optimal state transfer between localized modes by considering an initial Gaussian wave packet $\psi_0(x)$ which is transferred to a prescribed target $\psi_T(x)$, using Krotov’s method with different propagators. Two scenarios are considered here: (i) In the first scenario, the target is a symmetric state centered at the origin; a Gaussian wave packet with $x_0 = 0$ is taken as the initial state and is driven to the symmetric target using Krotov's method. (ii) In the second scenario, the target is an asymmetric state localized around the lattice site $x_0 = -25$; a Gaussian wave packet centered at $x_0 = -25$ is used as the initial state and is optimally transferred to this localized asymmetric target. The parametric values considered here correspond to a deep lattice with $V_0 = 10 E_R$ and a strong harmonic trap with strength given by $m \omega_T^2 d^2/2 = 0.00032E_R$, where $E_R\approx h \times 3.62~\mathrm{kHz}$ is the recoil energy of a single atom and $d=397.5nm$ is the lattice spacing.
	The optimization is performed with error tolerance $\varepsilon = 10^{-5}$ and with a maximal number of fifty Krotov iterations.  
	Three propagation schemes are compared:
	\begin{itemize}
		\item the fourth-order exponential commutator-free (CF-Exp) scheme \cite{AF:2011};
		\item the fourth-order Cayley-Magnus (M-Cayley) scheme \cite{Iserles:2001};
		\item the fourth-order commutator-free Cayley (CF-Cayley) method \cite{Wembe2025}.
	\end{itemize} 
	
	\begin{table*}[!ht]
		\centering
		\begin{tabular}{|l|l|l|l|r|r|}
			\hline
			Scheme & Initial state & Converged & CPU time (s) & Iterations & Final fidelity \\
			\hline 
			CF-Exp & Gaussian ($x_0=0$) & \cmark & 460.903 & 4 & 0.999981 \\
			Cayley–Magnus & Gaussian ($x_0=0$) & \cmark & 48.645 & 4 & 0.999981  \\
			CF–Cayley & Gaussian ($x_0=0$) & \cmark & 50.436 & 4 & 0.999981 \\
			\hline 
			CF-Exp  & Gaussian ($x_0=-25$) & \xmark & 5676.996 & 50 & 0.838997  \\
			Cayley–Magnus & Gaussian ($x_0=-25$) & \cmark & 122.791 & 10 & 0.734326  \\
			CF–Cayley & Gaussian ($x_0=-25$) & \cmark & 89.419 & 4 & 0.987149  \\
			\hline
		\end{tabular}
		\caption{Comparison of propagation schemes in the linear Schrödinger case with error \\ tolerance $\varepsilon = 10^{-5}$. A check mark (\cmark) indicates convergence of the Krotov iterations.}
		\label{tab:linear-results}
	\end{table*}

	\begin{figure}[!ht]
		\centering
		\def\sizefig{0.22}
		\includegraphics[width=\sizefig\textwidth]{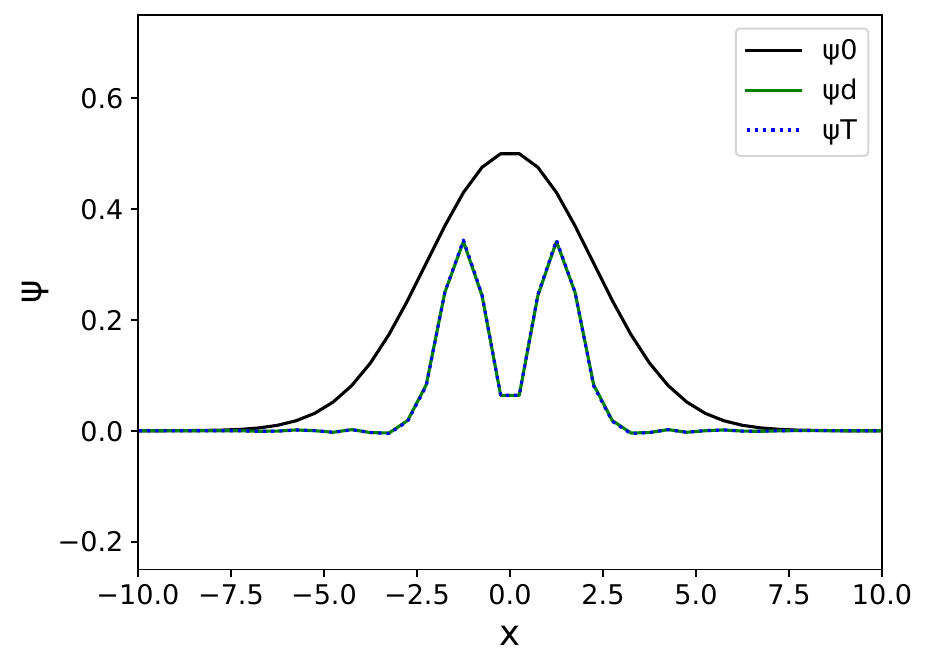}
		\includegraphics[width=\sizefig\textwidth]{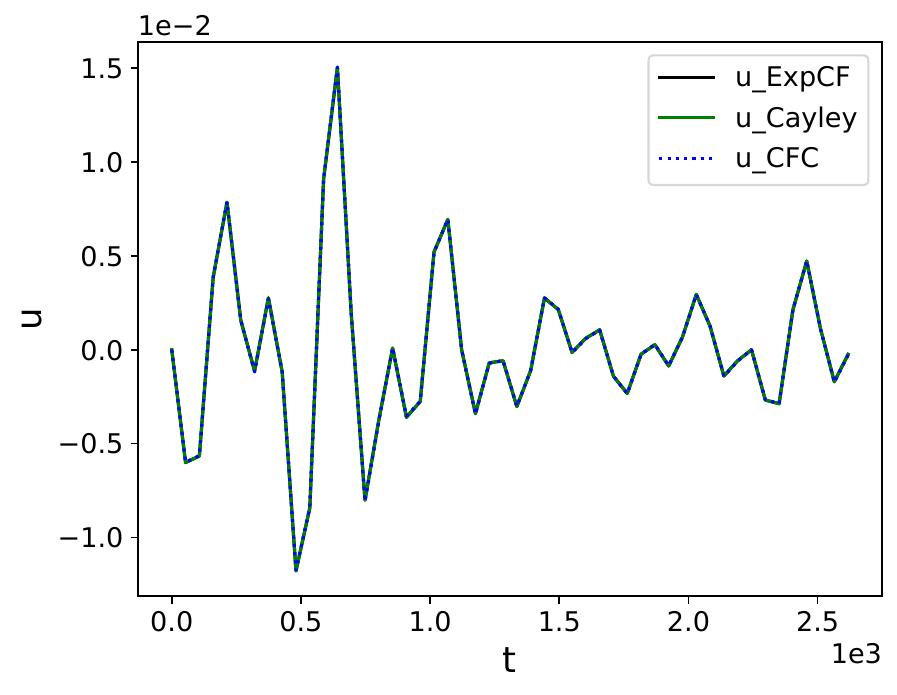}
		\includegraphics[width=\sizefig\textwidth]{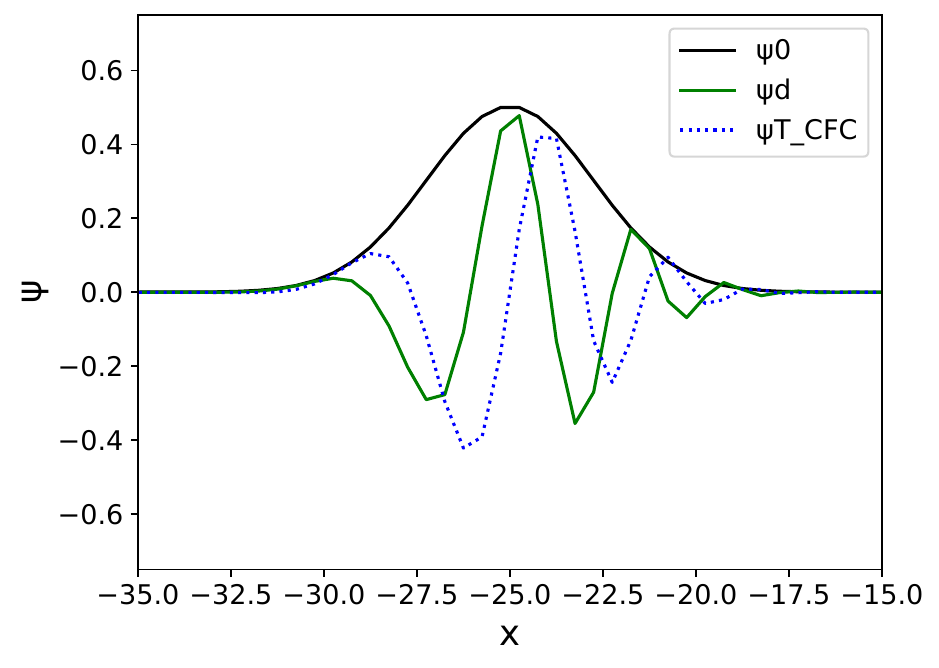}
		\includegraphics[width=\sizefig\textwidth]{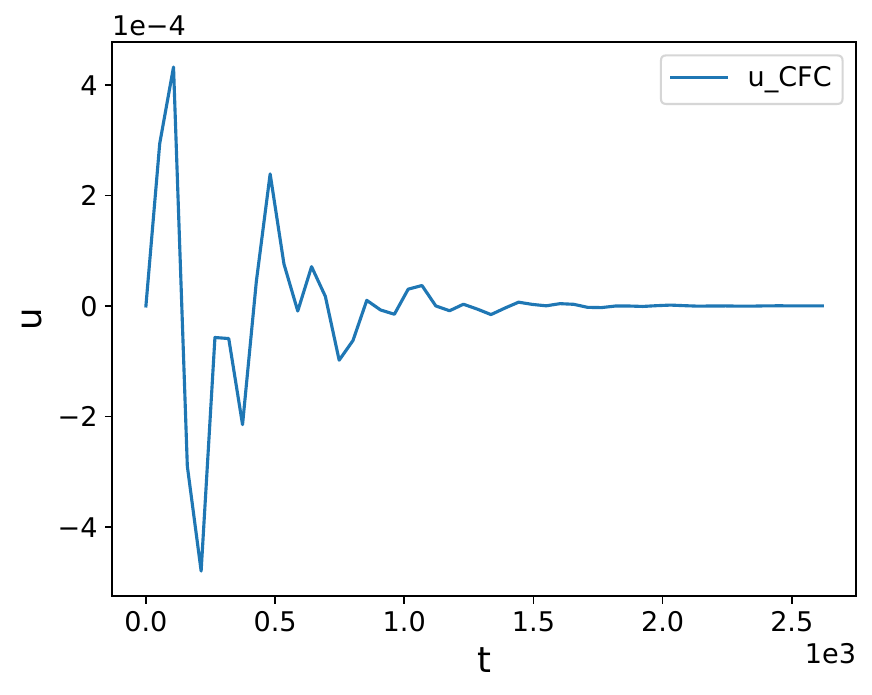}
		\caption{Quantum state transfer with the linear Schrödinger equation for non-interacting cold atoms in a driven parabolic lattice. The initial Gaussian wave packet at $x_0 = 0.0$, the target state, and the corresponding solutions of the Schrödinger equation are shown in the top-left panel, while the evolution of the control function is displayed in the top-right panel. The bottom-left panel presents a different target state obtained from an initial Gaussian wave packet centered at $x_0 = -25.0$, and the associated control function is shown in the bottom-right panel.
		}
		\label{fig:linear-case}
	\end{figure}
	
	\emph{Discussion.} Figure~\ref{fig:linear-case} presents the results for the optimal transfer in the two considered scenarios. The top-left panel depicts a high-fidelity transfer of the initial Gaussian wave packet at the origin to the symmetric target state, with the corresponding optimal control function shown in the top-right panel. In the bottom-left panel, the Gaussian wave packet at $x=-25$ retraces the overall shape of the asymmetric target state. The evolution of the control function in this case, shown in the bottom-right panel, exhibits a rapid modulation of the driving field. This modulation enables the wave packet to achieve the target shape in a shorter time, although some deviations from the desired state still remain.

	The performance of the three propagation schemes is summarized in Table~\ref{tab:linear-results}. For the symmetric target ($x_0 = 0$), all schemes converge to nearly identical final fidelity of about $0.999981$. However, the computational cost differs significantly: while the CF-Exp method requires the longest CPU time (about $460.9$~s), the Cayley–Magnus and CF–Cayley schemes achieve the same accuracy with CPU times of only $48.6$~s and $50.4$~s, respectively. For the asymmetric target ($x_0 = -25$), the performance difference becomes even more pronounced. The CF-Exp scheme fails to converge even after $50$ iterations, whereas both the Cayley–Magnus and CF–Cayley methods converge successfully. Among these, the CF–Cayley approach achieves the highest final fidelity ($0.987$) in only four iterations, with reduced computational time ($\sim89$~s) compared to the Cayley–Magnus method ($\sim123$~s, 10 iterations). These results demonstrate that the CFC (CF–Cayley) integrator preserves unitarity by construction, maintains high accuracy, and provides a speed-up of roughly an order of magnitude compared to the exponential-based method for the symmetric case, and over an order of magnitude for the asymmetric case, for comparable error tolerances, confirming its efficiency for repeated forward and backward propagation in Krotov’s algorithm. 
	
	\begin{remark}
		For the numerical examples considered here, no analytical solution of the Schrödinger equation is available. To ensure a fair comparison, the reference target state $\psi_d$ is computed by solving the Schrödinger equation with the fixed control $u(t)=\sin^2(t)$ using the same numerical propagator as that employed in the corresponding optimization method for the forward and backward propagations. This ensures that the reported fidelities in Table~\ref{tab:linear-results} reflect differences in the optimization performance rather than discrepancies arising from different propagation schemes.
	\end{remark}

	\subsection{Example 2: Nonlinear Schrödinger (Gross--Pitaevskii) Equation for Interacting Bose-Einstein Condensates}
	\label{sec:nonlinear-example}
	
	To assess the efficiency of the proposed integrators in a nonlinear setting, we consider the one-dimensional Gross--Pitaevskii equation (GPE)
	\small{
		\begin{eqnarray}
			i\partial_t \psi(x,t) =
			\left(-\Delta + V(x) + g|\psi(x,t)|^2 + u(t) W(x)\right)\psi(x,t),&& \nonumber \\
			\hspace{4.5cm} \psi(x,0)=\psi_0(x),&&
		\end{eqnarray} \label{eq:13}
	}
	where $V(x)$ represents a harmonic trapping potential (with $V(x) = x^4 - 10x^2$ in the numerical simulations, $u(t)W(x)$ an external control field (with $W(x) = 5x^2$ in numerical simulations), and $g \ge 0$ the nonlinear coupling constant describing mean-field interactions. 
	The term $g|\psi|^2$ makes the system nonlinear and introduces amplitude-dependent phase shifts, typical for weakly interacting Bose--Einstein condensates.  
	
	\begin{figure}[!ht]
		\centering
		\def\sizefig{0.47}
		\includegraphics[width=\sizefig\textwidth, height=0.2\textwidth]{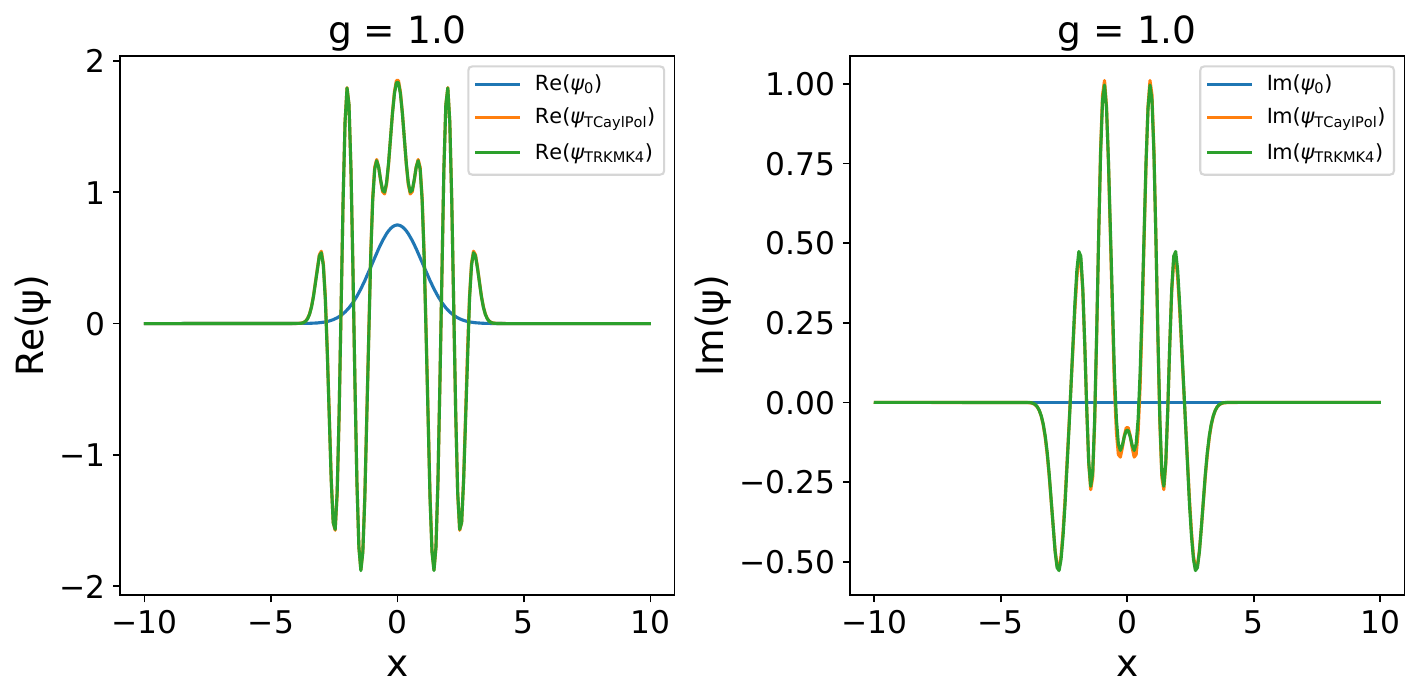}
		\caption{Quantum state transfer governed by the Gross-Pitaevskii equation under the fixed driving control $u(t)= u_c \sin(t)$. Starting from a Gaussian initial state localized at $x=0$, the system evolves toward the target state. The real and imaginary components of the target state are shown in the left-panel and right-panel, respectively, and are reproduced with high fidelity by both the Runge-Kutta-Munthe-Kaas (RKMK4) and Non-Linear CF-Cayley (CaylPol) methods.
		}
		\label{fig:non-linear-case}
	\end{figure}
	
	In this section, we fix the control field $u(t)= u_c \sin(t)$ and perform pure forward propagation using the CaylPol and the $4th$-order Runge-Kutta-Munthe-Kaas (RKMK4) methods. The optimal solution is shown in Fig.~\ref{fig:non-linear-case}. 
	Moreover we compare the CPU time required for one complete time integration for several values of the nonlinear parameter $g$, in Table~\ref{tab2:GPE-runtime}. For all tested values of $g$, the CaylPol method achieves almost the same accuracy as the RKMK4, while consistently requiring less computational time per integration step. The relative speed-up increases with $g$, since the nonlinear term amplifies stiffness in the dynamics and the unitary Cayley structure helps to stabilize the propagation.

	\begin{table}[!ht]
		\centering
		\begin{tabular}{|c|c|c|c|}
			\hline
			\textbf{Nonlinearity $g$} & 
			\textbf{CaylPol} & \textbf{RKMK4} \\
			\hline
			$0.0$   & $6.287$ & $343.798$  \\
			$0.5$   & $6.256$ & $430.991$  \\
			$1.0$   & $7.465$ & $343.148$  \\
			$2.0$   & $8.266$ & $490.429$  \\
			$10.0$  & $9.556$ & $580.881$  \\
			$20.0$  & $23.595$ & $658.118$ \\
			\hline
		\end{tabular}
		\caption{CPU-time comparison between RKMK4 and CaylPol schemes.}
		\label{tab2:GPE-runtime}
	\end{table}
	
	\begin{remark}
		(i) A detailed theoretical complexity analysis of the efficiency improvements observed in the numerical experiments is beyond the scope of this work. We refer to \cite{Wembe2025} for the theoretical analysis of the Cayley commutator--free integrator itself, including its order of accuracy and structural properties for quantum dynamical systems. In contrast to Magnus--based schemes and exponential commutator--free integrators, the proposed method avoids both commutator evaluations and matrix exponentials, replacing them by Cayley transforms and linear combinations of Hamiltonians, which reduces the computational cost per time step. A systematic complexity comparison of such structure--preserving integrators within optimal control algorithms will be the subject of future work.\\
		(ii) The purpose of Table~II is to compare the computational cost of the different propagation schemes for a fixed control function. Since the considered integrators are applied with the same discretization parameters, the resulting solutions are visually indistinguishable, and the table therefore focuses on the running time rather than on fidelity comparisons.
	\end{remark}

	\section{Conclusion}
	\label{sec:conclusion}
	
	We have introduced a family of structure-preserving numerical schemes for quantum optimal control based on commutator-free Cayley (CF-Cayley and CaylPol) integrators. Within the Krotov framework, these methods replace the usual exponential or Runge-Kutta propagators by unitary Cayley compositions that maintain 
	the geometric properties of quantum evolution at the discrete level. For linear systems, the CF-Cayley method achieves the same accuracy as 
	high-order exponential or Cayley-Magnus schemes while reducing computational cost.  
	In nonlinear settings such as the Gross--Pitaevskii equation, the approach preserves norm and stability even in strongly interacting regimes, providing a robust alternative to conventional explicit integrators.
	Beyond improved numerical efficiency, the CF-Cayley formulation offers a systematic way to design higher-order, structure-preserving propagators on Lie groups without commutator evaluations. Its simplicity and generality make it well suited for large-scale quantum optimal control problems, where repeated forward and backward propagations dominate computational cost. Future work will focus on adaptive step-size strategies, energy-preserving variants, and the inclusion of derivative-based cost functionals for constrained or amplitude-positive controls, as briefly discussed in Appendix~\ref{app:alt-cost}.

	\section*{APPENDIX}
	
	\subsection{Cayley-Magnus expansion method}
	\label{sec:Cayley-Magnus}
	The idea of the Cayley-Magnus expansion is to approximate the solution as the Cayley transform, 
	\begin{equation}
		\label{eq:Cayley_approx}
		\psi(t) = \Cay(\Omega(t)) = \left(I + \tfrac{1}{2} A(t)\right)^{-1} \left(I - \tfrac{1}{2} A(t)\right),
	\end{equation}
	of a series $\Omega(t) = \sum_{m=0}^\infty \Omega_m$. Given that $\psi$ is the solution of \eqref{eq:system_dynamics}, then $\Omega$ satisfies the following system.
	\begin{lemma}
		\label{lemma:1}
		Let $\psi(t)$ be the solution of system \eqref{eq:system_dynamics}, with $-1 \notin \sigma(\psi(t)\psi_0^{-1})$\footnote{where $\sigma(A)$ defined the spectrum of the matrix A.} for any $t \in [0,T]$. Let $\psi(t) = \Cay(\Om(t))\psi_0$, then $\Om$ is the solution of the system: 
		\begin{equation}
			\label{eq:skew-Herm-prob}
			\dot{\Om}(t) = A(t) - \frac{1}{2}[\Om, A(t)] - \frac{1}{4}\Om A(t) \Om, \quad \Om(t_0) = \Om_0, 
		\end{equation}
		with the Lie bracket (commutator) of two matrices $A$ and $B$ defined by $[A,B] = A\cdot B - B\cdot A$.
	\end{lemma}
	By substituting 
	$\Omega(t) = \sum_{m=0}^\infty \Omega_m$ into \eqref{eq:skew-Herm-prob} and integrating over $[0,T]$ one can recursively computed the $\Om_j$ (see \cite[Sec.~2.2]{Wembe2025} for more details) as follows:
	\begin{eqnarray}
		\Om_1(t) &=& \int_0^t A(\xi) \diff \xi, \nonumber \\
		\Om_2(t) &=& \hspace{-0.2cm} - \frac{1}{2} \int_0^t [\Omega_1(\xi), A(\xi)] \diff \xi,\\
		\Om_m(t) &=& ~ \frac{1}{2} \int_0^t [\Om_{m-1}(\xi), A(\xi)] \diff \xi  \nonumber \\ 
		&& - \frac{1}{4} \sum_{k=1}^{m-2} \int_0^t \Om_{m-k-1}(\xi)  A(\xi) \Om_k(\xi) \diff \xi, ~ \text{for} ~ m \geq 3. \nonumber
	\end{eqnarray}

	\begin{lemma}
		Truncating the Cayley-Magnus expansion $\Omega(t) = \sum_{m=0}^\infty \Omega_m$ at a given order $p$, i.e.~setting
		\begin{equation}
			\label{eq:Cayley-truncate-expension}
			\Om(t) \approx \Om_p(t) = \sum_{m=1}^{p} \Om_p,
		\end{equation}
		leads to a $p$-order approximation.
	\end{lemma}
	Combining \eqref{eq:Cayley_approx} and \eqref{eq:Cayley-truncate-expension} one can build any higher-order Cayley-Magnus expansion scheme. Moreover, the order in \eqref{eq:Cayley-truncate-expension} is
	increased to $p+1$ when we consider specifically the Gauss-Legendre quadrature to approximate the integrals in $\Om_j$. One has for instance the following $4th$-order Cayley-Magnus expansion scheme (CME), see \cite{Iserles:2001}.
	
	\begin{lemma}
		\label{lemma:2}
		Given a time step $\delta t$, the following approximation holds,
		\begin{equation}
			\label{eq:simplify-approx}
			\begin{aligned}
				Y(\delta t) &= \Cay(\Om(\delta t)) \\ 
				&= \Cay\left(A_1 -\frac{1}{6} [A_1, A_2] - \frac{1}{12} A_1^3 \right)  + O(\delta t^5), 
			\end{aligned}
		\end{equation}
		with
		\begin{equation}
			\label{eq:Gauss-Legendre}
			\begin{aligned}
				A_1 &= \frac{\delta t}{2}(A^1 + A^2) + O(\delta t^5), \\
				A_2 &= \frac{\delta t\sqrt{3}}{2}(A^2 - A^1) + O(\delta t^5),
			\end{aligned}
		\end{equation}
		and 
		\[
		\begin{aligned}
			A^1 &= A\left(\left(\frac{1}{2}-\frac{\sqrt{3}}{6} \right) \delta t\right), \\
			A^2 &= A\left(\left(\frac{1}{2}+\frac{\sqrt{3}}{6} \right) \delta t \right).
		\end{aligned}
		\]
	\end{lemma}

	\subsection{Commutator-Free Cayley (CFC) method}
	\label{sec:CFC-method}
	
	The CFC idea is to approximate the propagator by a product of Cayley transforms,
	\begin{equation}
		\label{eq:Comm_free}
		Y(\delta t) \approx \prod_{\ell=1}^s \Cay(\widetilde\Om_\ell(\delta t)),
	\end{equation}
	where the effective operators $\widetilde{\Om}_\ell = \sum_{k=1}^2 \alpha_{\ell,k} A_k$ are linear combinations of $A$ evaluated at quadrature nodes in $[t,t+\delta t]$ and $\alpha_{i,k} \in \R$ are some coefficients, then use the Cayley-BCH formula to determined, when there exists, such that approximation \eqref{eq:Comm_free} leads to a desired higher-order approximation.  

	As an example, a symmetric three-stage ($s=3$) fourth-order CFC scheme reads
	\begin{equation}
		U^{[4]}_{\mathrm{CFC}}(\delta t) 
		= \Cay(\delta t, \widetilde{A}_1)\,
		\Cay(\delta t, \widetilde{A}_2)\,
		\Cay(\delta t, \widetilde{A}_3),
	\end{equation}
	with
	\begin{align}
		\widetilde{A}_1 &= \alpha_{11} A(t+c_1 \delta t) 
		+ \alpha_{12} A(t+c_2 \delta t), \nonumber \\
		\widetilde{A}_2 &= \alpha_{21} A(t+c_1 \delta t) 
		+ \alpha_{22} A(t+c_2 \delta t), \\
		\widetilde{A}_3 &= \alpha_{31} A(t+c_1 \delta t) 
		+ \alpha_{32} A(t+c_2 \delta t), \nonumber
	\end{align}
	where the nodes and coefficients are
	{\small
		\[
		\begin{cases}
			c_{1,2} = \tfrac{1}{2} \mp \tfrac{\sqrt{3}}{6}, 
			~~ \alpha_{22} = 0, \\[0.3em]
			\alpha_{11} = \dfrac{2^{1/3}}{3} + \dfrac{2^{2/3}}{6} + \dfrac{2}{3},
		\end{cases}
		\begin{cases}
			\alpha_{12} = - \alpha_{32} = \alpha_{11} - \alpha_{11}^2, \\
			\alpha_{31} = \alpha_{11}, ~~ 
			\alpha_{21} = 1 - 2\alpha_{11}.
		\end{cases}
		\]
	}
	%

	\subsection{Discussion on Alternative Cost Functions}
	\label{app:alt-cost}
	
	In several works on quantum and mean-field optimal control 
	(see e.g.~\cite{Jager2014, Reich2012, Grond2009}), 
	the cost functional is augmented by an additional penalty term involving the derivative of the control,
	\begin{equation}
		J(\psi(T),u) 
		= J_T(\psi(T)) + \frac{\alpha}{2}\int_0^T |\dot{u}(t)|^2\,\mathrm{d}t,
	\end{equation}
	where $J_T(\psi(T)) = \frac{1}{2}\left(1 - |\langle \psi_T, \psi(T)\rangle|^2 \right)$ measures the terminal infidelity and $\alpha>0$ balances control smoothness against target accuracy. This term penalizes rapid temporal variations of the control and is commonly used in GRAPE-type schemes~\cite{Grond2009, Reich2012}.  
	
	From a physical perspective, such regularization improves experimental feasibility 
	by enforcing smooth and bounded control fields.  
	Mathematically, it ensures the well-posedness of the optimization problem 
	in $H^1(0,T)$ rather than in $L^2(0,T)$, 
	and leads to an additional Euler–Lagrange equation for the optimal control, 
	typically of the form
	\[
	-\ddot{u}(t) 
	= \frac{1}{\alpha}\mathrm{Im}\langle \lambda(t),\,\partial_{u}H(u(t))\,\psi(t)\rangle,
	\]
	which couples the control curvature to the forward and adjoint dynamics.  
	An additional advantage of this formulation is that, 
	by optimizing over the time derivative $\dot{u}(t)$ 
	and reconstructing $u(t)$ through integration, 
	the positivity of the control field can be naturally enforced. This approach is particularly relevant when the control represents a physical intensity or a laser amplitude that must remain nonnegative. Although such derivative-based regularization has not been explored in the present work, it provides an interesting direction for future developments of structure-preserving optimal control schemes.

	\section*{ACKNOWLEDGMENT}
	
	The authors would like to thank P.~Singh and C.~Offen for many valuable discussions and insights on the theoretical aspects of this work. We are also grateful to F.~Geiesler for assistance with parts of the numerical implementation and code development.


\end{document}